\documentclass[article]{amsart}
\usepackage{graphicx}
\newtheorem{Theorem}{Theorem}[section]

\newtheorem{Proposition}[Theorem]{Proposition}

\theoremstyle{definition}

\title{
A double large deviation principle for Monge-Amp\`ere gravitation
} 
\author{Yann Brenier, 
\\
CNRS UMR 7640, 
\\
Ecole Polytechnique, Palaiseau, France.
} 
\date\today

\begin{document}

\begin{abstract}
\noindent
Monge-Amp\`ere gravitation is a nonlinear modification of classical Newtonian gravitation,
when the Monge-Amp\`ere equation substitutes for the Poisson equation.
We establish, through two applications of the large deviation principle,
that the MA gravitation for a finite number of particles can be reduced, through a double
application of the large deviation principle, to the simplest possible
stochastic model: a collection of independent Brownian motions with vanishing noise. 
\end{abstract}
\maketitle

\subsubsection*{This paper is dedicated to Professor Tai-Ping Liu for his seventieth birthday}

\section*{Introduction}
The purpose of this paper is twofold. We first want to make a short presentation of the
Monge-Amp\`ere gravitational (MAG) model, which has been introduced in \cite{Br8},
in close connection with earlier works such as \cite{Br4,BL,Lo,CGP,AG} as well
as with optimal transport theory (see \cite{Vi,AGS}). 
This model can be seen as a nonlinear modification of the classical model
of Newtonian gravitation, for which we use the fully nonlinear Monge-Amp\`ere equation 
as a substitute for
the linear Poisson equation to derive the gravitational potential from the density field. We also
briefly compare the Monge-Amp\`ere and the Newton gravitational models, emphasizing
that in the MAG model there is an absolute control of the acceleration of the gravitating
particles, no matter how concentrated they can be, in sharp contrast with
Newtonian gravitation, for which some particles may runaway to infinity in finite time.
However, the main purpose of this paper is about deriving the MAG model from 
a completely elementary microscopic model in which a finite number of particles just move as 
independent Brownian trajectories without any interaction. In order to get the MAG model,
we need two applications of the large deviation principle (LDP) \cite{FW,DZ}.
(Let us mention the interesting connection between 
large deviation principles, gradient flows and
optimal transport theory recently made in \cite{Le,MPR}, 
which was certainly influential for us.)
Through the first application of the LDP, we get 
a first order (in time) dynamical system based on the indistinguishability of the gravitating
particles. (In a rather paradoxical way, the indistinguishability principle leads
to a model where particles do interact.) The second application of the LDP enables us to
lift this first order dynamical system to a second order one which is nothing but the discrete
version of the MAG model.
Let us emphasize that this derivation is purely formal and further investigations are clearly
needed to get a complete and rigorous theory. 
The first section of this paper will be devoted to a presentation
of both the classical Newtonian and the Monge-Amp\`ere models, with a suitable formulation
of the MAG model relying on "optimal transport" tools \cite{Br1,Br2,Vi}. 
The second section is devoted to the
double application of the LDP in order to recover the discrete version of the MAG model.

\section{\bf{Monge-Amp\`ere and Newton gravitations}}
\subsection{Classical Newtonian gravitation}
To describe classical Newton gravitation, 
let us attach to each gravitating particle a label $a$,
which we suppose to belong to an abstract Borel measure space $(\mathcal{A},\lambda)$,
and its position at time $t$, $X_t(a)\in \mathbb{R}^d$ (classically $d=3$). 
For each time $t$, we denote by $\rho_t$ the image measure of $\lambda$ by $X_t$,
defined on $\mathbb{R}^d$ by
\begin{equation}\label{image}
\int_{\mathbb{R}^d} \gamma(x)\rho_t(dx)=\int_{\mathcal{A}} \gamma(X_t(a))\lambda(da)
\;\;\;\forall \gamma\in C_c^0(\mathbb{R}^d).
\end{equation}
For the sake of simplicity, we further assume $\rho_t(dx)$ to be $\mathbb{Z}^d$ periodic in $x$ and of unit mean
on the unit cube $[0,1]^d$, so that we can see $\rho_t(dx)$ as a probability measure
on the flat torus 
$\mathbb{T}^d=\mathbb{R}^d/\mathbb{Z}^d$.
(This kind of assumption is common in computational Cosmology \cite{FMMS}.)
We are now ready to write the Newtonian model:
\begin{equation}\label{newton}
\frac{d^2 X_t(a)}{dt^2}=-\nabla \phi_t(X_t(a)),
\;\;\;\bigtriangleup \phi_t=\rho_t-1,
\end{equation}
where the average density $1$ has been subtracted out from the field $\rho_t$ so 
that the "gravitation potential" $\phi_t$ is well defined from the Poisson equation as a zero mean
$\mathbb{Z}^d-$periodic function. (This is a perfectly meaningful assumption: there should be no resulting force $\nabla\phi_t$ out of a uniform $\rho_t$ .)
If, at time $t$, $\rho_t$ is just a plain probability measure on $\mathbb{T}^d$,
then $\nabla\phi_t$, obtained through the Poisson equation $\bigtriangleup \phi_t=\rho_t-1$, is merely an $L^p$ function for $p<(1-d^{-1})^{-1}$ on $\mathbb{T}^d$. This makes very dubious the meaning of
ODE $$\frac{d^2 X_t(a)}{dt^2}=-\nabla \phi_t(X_t(a))$$ in the general case.
This is why Newtonian's gravitation usually splits up into two, rather unrelated theories,
according to the choice of its initial conditions:
\\
1) As $\mathcal{A}$ is a discrete set of $N$ particles,
then Newtonian gravitation boils down to an "$N$-body problem", set on the flat torus $\mathbb{T}^d$,
which is of form
\begin{equation}\label{N-body}
\frac{d^2 X_t(a)}{dt^2}=R(X_t(a))+\sum_{b\ne a}G(X_t(a),X_t(b)),
\end{equation}
where $R$ is a fixed function depending on the torus, 
and the singular kernel $G(x,y)$, associated with the operator 
$-\nabla\bigtriangleup^{-1}$,
behaves as $(y-x)/|y-x|^{d}$, as $x$ approaches $y$. 
This "bad" singularity makes the study of the "$N$-body problem" very hard, in particular because of the possibility that particles may collide or runaway at infinity, in finite time. (See \cite{Ch} as a recent reference.)
\\
2) As $\mathcal{A}$ corresponds to a "continuum" of particles, it is fruitful
to introduce, at each time $t$, the nonnegative measure ${f}_t$ defined on the
"phase space" $\mathbb{T}^d\times\mathbb{R}^d$ by
\begin{equation}\label{kinetic}
\int_{\mathbb{T}^d\times\mathbb{R}^d
} \gamma(x,v)f_t(dx,dv)
=\int_{\mathcal{A}} \gamma(X_t(a),\frac{dX_t(a)}{dt})\lambda(da),
\;\;\;\forall \gamma\in C^0_c(\mathbb{T}^d\times\mathbb{R}^d).
\end{equation}
A straightforward calculation shows that $f$ must satisfy, at least formally,
the "Vlasov-Poisson" system \cite{Re}
\begin{equation}\label{VP}
\partial_t f_t+\nabla_x\cdot(vf_t)-\nabla_v\cdot(\nabla\phi_tf_t)=0,
\;\;\;\bigtriangleup \phi_t=\rho_t-1=\int_{\mathbb{R}^d}f_t (dv)-1.
\end{equation}
The meaning of this equation is unclear in the general case,
because of the dubious product of $\nabla\phi_t$ by $f_t$.
However, it is not hard to see that
$0\le f_t(dx,dv)\le Cdxdv$ for some positive
constant $C$ is a consistent property for equation (\ref{VP}). It turns out that this property is sufficient
to justify the multiplication of $f_t$ by $\nabla\phi_t$ at least when the initial condition $f_0$ vanishes
for large values of $v$, say for $|v|\ge C'$ for some other constant $C'$. Then, the Newtonian
model gets well defined and, at least for $d\le 3$, global weak solutions can be shown to exist 
globally in time, as $f_0$
is given, whatever are constants $C$ and $C'$. In addition, such solutions are unique and
smooth as long as $f_0$ has a smooth density with respect to the Lebesgue measure.
(See \cite{Pf,LP,BGP} for such results, stated on the whole space $\mathbb{R}^d$
rather than on $\mathbb{T}^d$. See also \cite{Re}, as well as
\cite{DLYY} for the closely related
Euler-Poisson system.)

\subsection{Monge-Amp\`ere gravitation}

We now introduce the Monge-Amp\`ere gravitation (MAG) model which differs from (\ref{newton}) 
just by the substitution of the
Monge-Amp\`ere equation for the Poisson equation:
\begin{equation}\label{MAG}
\frac{d^2 X_t(a)}{dt^2}=-\nabla \phi_t(X_t(a)),
\;\;\;\det(\mathbb{I}+D^2\phi_t)=\rho_t,
\end{equation}
where $\mathbb{I}$ denotes the $d\times d$ identity matrix and $\phi_t$ is
solution to
the Monge-Amp\`ere equation (in a suitable sense), which is
 $\mathbb{Z}^d-$periodic with zero mean and satisfies
the (weak) ellipticity condition $\mathbb{I}+D^2\phi_t(x)\ge 0$,
in the sense of symmetric matrices, for every $x$. [Notice that $\phi_t$ is unique and Lipschitz
continuous (resp. smooth) as soon as $\rho_t(dx)$ has an integrable (resp. smooth and positive) 
density with respect to the Lebesgue measure $dx$, see \cite{Co} for example.]
We see that Newtonian gravitation can be formally recovered from the MAG model just by expanding the determinant about $\mathbb{I}$ and retaining only the linear part:
$$
\det(\mathbb{I}+D^2\phi_t)\sim 1+{\rm{Trace}}D^2\phi_t=1+\bigtriangleup \phi_t.
$$
Notice that, as $d=1$ (which is a case of limited interest, describing 
gravitating ``parallel pancakes''), the MAG model coincide with Newtonian
gravitation. (As a consequence, our derivation of the MAG model from
a double application of the large deviation principle, obtained in the
second section of the present paper, is also valid
for the Newtonian gravitational model in one space dimension.) However,
let us emphasize that,
to the best of our knowledge, MAG has never been considered by any physicist
in dimension larger than one.
It has been so far only considered by mathematicians (see \cite{Br8}
and \cite{BL,Lo,CGP,AG} for closely related topics), mostly because of its close connection
with optimal transport theory (as discussed in the next subsection).
It is easy to describe, as we did for the Newtonian model,
the MAG model through a "kinetic equation", the so-called "Vlasov-Monge-Amp\`ere" (VMA) system 
\cite{BL,Lo}:
\begin{equation}\label{VMA}
\partial_t f_t+\nabla_x\cdot(vf_t)-\nabla_v\cdot(\nabla\phi_tf_t)=0,
\;\;\;\det(\mathbb{I}+D^2\phi_t)=\rho_t
=\int_{\mathbb{R}^d}f_t (dv).
\end{equation}
As for the Vlasov-Poisson system (\ref{VP}), the existence of global weak solutions can be
shown as soon as 
$$
0\le f_0(dx,dv)\le Cdxdv,\;\;\;f_0(\mathbb{T}^d\times \{v\in \mathbb{R}^d,\;\;|v|\le C'\})=0,
$$
for some positive constants $C,C'$. 
Existence of a unique smooth solution, but only for a short time interval, 
has been proven by Loeper \cite{Lo} provided $f_0(dx,dv)$ 
(resp. $\rho_0(dx)$) has a smooth density with respect to $dxdv$ (resp. $dx$) and is
uniformly compactly supported in $v\in\mathbb{R}^d$ (resp. is strictly positive).
These are still limited results, due to the full non-linearity of the Monge-Amp\`ere
equation that leads to analytic difficulties, in particular when seeking for smooth solutions.
However, the Monge-Amp\`ere equation enjoys remarkable properties, closely related to the
theory of optimal transportation \cite{Br2,Vi}. This is why we are going to introduce a related
formulation of the MAG model, with interesting geometric features.

\subsection{The MAG model written in optimal transportation terms}
\label{subsection-MAG}

For the description of the MAG model in optimal transport terms, it is convenient to discuss
the model in a (slightly) different and more abstract framework. We consider a metric measure
space $(\mathcal{A},\lambda)$, made of a compact subset $\mathcal{A}$ of $\mathbb{R}^d$
equipped with a Borel probability measure $\lambda$. Two typical examples
are, on one side, the unit cube with the Lebesgue measure, and, on the other side,
any set of $N$ points equipped with the (normalized) counting measure. In the first case,
we will speak of the "continuous" case, while, in the second case, we will speak of the
"discrete" case. (Of course many others situations could be also considered, in particular
the flat torus $\mathbb{T}^d$ as we did in the previous subsections, but we will
focus on these two cases only.)
We introduce
the separable Hilbert space $H$ 
of all $\lambda-$square-integrable maps from $\mathcal{A}$ to $\mathbb{R}^d$,
$H=L^2(\mathcal{A},\lambda;\mathbb{R}^d),$
with norm and inner product respectively denoted by $||\cdot||$ and $((\cdot,\cdot))$.
(Notice that,
when $\mathcal{A}$ is a finite
set of $N$ points in $\mathbb{R}^d$,  then $H\sim \mathbb{R}^{Nd}$ is of finite dimension.)
We crucially consider the subset $S$ of
all  measure-preserving maps $s$ of $A$:
\begin{equation}
\label{VPM}
S=\{s\in H,\;\;\int_{\mathcal{A}} \gamma(s(a))\lambda(da)=
\int_{\mathcal{A}} \gamma(a)\lambda(da),\;\;\forall \gamma\in C^0(\mathbb{R}^d)\}.
\end{equation}
(Observe, in the discrete case, when $\mathcal{A}$ is made of $N$ distinct points 
$A(a)\in\mathbb{R}^d$, for $a=1,\cdots,N$,
$S$ just reads $S=\{(A(\sigma(1)),\cdots,A(\sigma(N))\in H,\;\;\sigma\in\mathcal{S}_N\}$
with $N!$ elements,
where $\mathcal{S}_N$ denotes the group of all permutations of the first $N$ integers.)
\\
\\
According to Edelstein's theorem \cite{Ed,Au}, 
given a separable Hilbert space $H$ and a closed bounded subset $S$,
there is, in the  sense of Baire, a generic set (i.e. containing a countable intersection of dense open subsets of $H$) of points $X$ for which 
there exists a unique closest point $\pi(X)$ on $S$. In
addition $\pi(X)$ is nothing but the gradient, at point $X$,
of the Lipschitz
convex function $\Pi$ defined on the Hilbert space $H$ by
\begin{equation}
\label{Pi}
\Pi(X)=\sup_{s\in S}((X,s))-\frac{||s||^2}{2},\;\;\;
\pi(X)=\nabla\Pi(X)={\rm{Arg}}\inf_{s\in S}\frac{||X-s||^2}{2}.
\end{equation}
In our particular case, we can say much more in the continuous case, when $\mathcal{A}$ is the unit cube 
with $\lambda$ as the Lebesgue measure, thanks to "optimal transport theory"
\begin{Theorem}(\cite{Br2})
\label{theo-polar}
Let $X\in H$ be a non degenerate map, in the sense that the image measure $\rho$
of $\lambda$ by $X$ is absolutely continuous with respect to the Lebesgue measure
on $\mathbb{R}^d$.
Then, $X$ has a unique closest point $\pi(X)$ on $S$. We also have
\begin{equation}
\label{defpi}
\pi(X)=T\circ  X,
\end{equation}
where $T$ is 
the unique map (in the $\rho$ a.e. sense) 
$y\in\mathbb{R}^d\rightarrow T(y)\in\mathcal{A}$ such that:
\\
i) there is a convex Lipschitz function 
$\psi:\mathbb{R}^d\rightarrow \mathbb{R}$, $\rho-$a.e.
differentiable, with:
\begin{equation}
\label{weakMA1}
T(y)=\nabla\psi(y),\;\;\;\rho-a.e.\;y\in\mathbb{R}^d
\end{equation}
ii) $\lambda$ is the image of $\rho$ by $T$:
\begin{equation}
\label{weakMA2}
\int_\mathcal{\mathbb{R}^d} \gamma(T(y))\rho(dy)
=\int_\mathcal{A} \gamma(a)\lambda(da),
\;\;\;\forall \gamma\in C^0(\mathbb{R}^d).
\end{equation}
\end{Theorem}

As a matter of fact, equations (\ref{weakMA1},\ref{weakMA2}) 
form a generalized formulation of the Monge-Amp\`ere problem on 
$\mathbb{R}^d$
\begin{equation}\label{MAP}
\rho =\det(D^2_x\psi),\;\;(\nabla\psi)({\rm{support}}(\rho))=\mathcal{A},
\;\;D^2\psi\ge 0,
\end{equation}
with a unique solution $\nabla\psi$ (in the $\rho-$a.e. sense).
(See \cite{Br1,Br2,Vi} for more details.)
This generalized formulation of the Monge-Amp\`ere equation
allows us to write the MAG model 
(\ref{MAG}) in a much more geometric way. Indeed, (\ref{MAG}) just reads
$$
\frac{d^2 X_t(a)}{dt^2}=X_t(a)-\nabla \psi_t(X_t(a)),
\;\;\det(D^2\psi_t)=\rho_t,
$$
after setting $\psi_t(x)=|x|^2/2+\phi_t(x)$, that we complete 
with the weak
ellipticity condition $D^2\psi_t\ge 0$ and the range
condition $\nabla\psi_t({\rm{support}}(\rho_t))=\mathcal{A}$ (which substitutes for
the $\mathbb{Z}^d$-periodicity condition we used in writing (\ref{MAG})).
Assume that, at time $t$, $X_t$ is non degenerate (or, in other words, $\rho_t$
is absolutely continuous 
with respect to the Lebesgue measure on $\mathbb{R}^d$). Then, using Theorem
\ref{theo-polar}, we may write
$$
\nabla\psi_t\circ X_t=\pi(X_t)=\nabla\Pi(X_t)
$$
and finally obtain
\begin{equation}
\label{MAG2}
\frac{d^2 X_t}{dt^2}=X_t-\pi\circ X_t=X_t-\nabla\Pi\circ X_t,
\end{equation}
where $\Pi$ is the Lipschitz convex function defined by (\ref{Pi}), $\pi$ its gradient and $S$ is the
set of measure preserving maps defined by (\ref{VPM}).
In the rest of this paper, we will retain (\ref{VPM},\ref{Pi},\ref{MAG2}) as our definition
of the MAG model.

\subsection{The discrete Monge-Amp\`ere gravitational model}

The geometric formulation (\ref{VPM},\ref{Pi},\ref{MAG2}) of the MAG
model is very convenient to get its discrete version,
as $\mathcal{A}$
is just a subset of $N$ points 
$A(a)$ in $\mathbb{R}^d$, for $a=1,\cdots,N$,
equipped with the counting measure,
in which case, we will speak of the "discrete MAG model with $N$ particles".
Indeed, at the discrete level, 
a time-dependent map $X_t$ can be seen as a sequence of $N$ "particles",
with positions $X_t(a)$, moving in $\mathbb{R}^d$,  for $a=1,\cdots,N$. (Notice
that $X_t$ gets degenerate at time $t$ in case of collisions.)
In this discrete setting, the MAG model
(\ref{VPM},\ref{Pi},\ref{MAG2}) reads as the finite dimensional dynamical system:
\begin{equation}
\label{MAG3}
\frac{d^2 X_t(a)}{dt^2}=X_t(a)-A(\sigma_t(a)),
\;\;
\sigma_t={\rm{Arg}}\inf_{\sigma\in\mathcal{S}_N}
\sum_{a=1}^N |X_t(a)-A(\sigma(a))|^2,
\end{equation}
where $\sigma\in\mathcal{S}_N$ denotes the set of all permutations of the first $N$
integers.
So we see that, at the discrete level, the MAG looks both very simple and very different from the classical
Newtonian $N$ body problem (\ref{N-body})! (Similar systems have been previously studied in
\cite{Br4,CGP}.)
\\
\\
The MAG model, as defined by (\ref{VPM},\ref{Pi},\ref{MAG2}), can be seen
as a dynamical system in a Hilbert space $H$ with a force term $\mathcal{F}(x)=x-\nabla\Pi(x)$
with a trivial linear part and the gradient of a Lipschitz convex function $\Pi$.
To the best of our knowledge, there is no theory for such an ODE in infinite dimension.
Indeed, we are very far from the standard Cauchy-Lipschitz setting.
This is a challenging open problem.
(See the related theory developed
by Ambrosio and Gangbo for some infinite dimensional hamiltonian systems
\cite{AG}. As a matter of fact, their main example is very similar to the MAG model, written in
a different way. See also \cite{BL,GNT}.)
However, for the discrete MAG model with $N$ particles, $H\sim \mathbb{R}^{Nd}$
is of finite dimension. Thus, equation (\ref{MAG2}), which also reads (\ref{MAG3}),
can be neatly solved
globally in time in the sense of Bouchut and Ambrosio \cite{Bo,Am}, 
for every initial condition $(X_0,\frac{dX_0}{dt})$, except on a negligible subset of
the  "phase space" $H\times H$, of zero $2Nd$ dimensional
Lebesgue measure. (A more accurate statement can be found in \cite{Am}.)
Indeed, the force term $\mathcal{F}(x)=x-\nabla\Pi(x)$  in the right-hand side of (\ref{MAG2}) is 
a smooth perturbation of a "bounded variation" function, since $\Pi$ is 
a Lipschitz convex function (which implies that $D^2\Pi$ can be seen as
a bounded nonnegative measure valued in the convex cone of nonnegative symmetric matrices,
and, therefore, that $\nabla\Pi$ is of bounded variation).
Notice that the exceptional set of bad initial conditions in the phase space is not empty, as it follows clearly from formulation
(\ref{MAG3}), where we see that the evolution of the particles becomes
ambiguous as different particles depart from the same position with $exactly$ the
same velocity.

\subsection{Monge-Amp\`ere versus Newtonian gravitational models}

It is now interesting to compare the Monge-Amp\`ere and the Newton gravitational models.
We first observe that, according to the MAG model, 
particles may never runaway to infinity in finite time. 
Indeed, from the optimal transportation formulation (\ref{VPM},\ref{Pi},\ref{MAG2}),
we immediately get 
$$
|\frac{d^2 X_t}{dt^2}-X_t|\;\le\; R=\sup_{a\in \mathcal{A}}\;|a|
$$
and deduce that 
$$
|X_t|+|\frac{dX_t}{dt}|\le 
(|X_0|+|\frac{dX_0}{dt}|+1)C\cosh t,\;\;
\forall t\in\mathbb{R},
$$
where $C$ depends only on $R$. 
Of course, there is nothing similar with Newtonian gravitation. Indeed, 
the Poisson equation in (\ref{newton}) does not behave well when particles concentrate.
For instance, if some particles concentrate as a delta measure at some point $y$ at time $t$,
then $\nabla\phi_t(x)$ has a singularity as bad as $(x-y)/|x-y|^{d}$.
\\
In this way, the MAG model has a lot of similarity with the Born-Infeld (BI) theory of the electromagnetic
field \cite{BI,BDLL,Br7,BY}, in which any electrostatic field is bounded by a universal constant.
(Notice that the BI model, which goes back to 1934, is no longer for use in Electrodynamics, but has enjoyed a remarkable revival in String Theory since the 1990s \cite{Po}.)
\\
In addition, the MAG model enjoys good properties of approximations by finite sums of
Dirac measures (see \cite{CGP} for closely related results). 
In sharp contrast, such discrete approximations have never been justified, to the best of our knowledge, in the case of 3D Newtonian gravitation, because of the bad singularities of the Green function for the Poisson equation \cite{HJ}. As a matter of fact, the treatment of point particles in classical Electrodynamics remains an outstanding open problem in both Theoretical and Mathematical Physics (see \cite{Fe,Sp} for instance).

\section{\bf{Formal derivation of the discrete Monge-Amp\`ere gravitation model from
a double application of the large deviation principle}}
 
In this second section, we provide a $formal$ derivation of the discrete
MAG model (\ref{MAG3}) from a double application
of the large deviation principle \cite{FW,DZ}.

\subsection{A basic model of $N$ independent Brownian particles} 
We fix a positive integer $N$ and a finite set $\mathcal{A}$ of 
$N$ distinct points $A(a)\in \mathbb{R}^d$, for $a=1,\cdots,N$.
Using the notations of subsection \ref{subsection-MAG},
we set $H=\mathbb{R}^{Nd}$ and 
\begin{equation}
\label{S2}
S=\{(A(\sigma(1)),\cdots,A(\sigma(N))\in H,\;\;\sigma\in\mathcal{S}_N\},
\end{equation}
where $\mathcal{S}_N$ denotes the set
of all possible permutations of the first $N$ integers.
\\
Given $\sigma_0\in \mathcal{S}_N$, 
we consider the motion of $N$ "particles" that move in the Euclidean space
$\mathbb{R}^d$, according to
\begin{equation}\label{brownian}
X^\varepsilon_t(a)=s_0(a)
+\varepsilon B_t(a),
\;\;s_0(a)=A(\sigma_0(a)),
\;\;\forall a=1,\cdots,N.
\end{equation}
where
$
(t\in \mathbb{R}_+\rightarrow B_t(a)\in\mathbb{R}^d)_{a=1,\cdots,N}
$ 
denote $N$ independent
realizations of the "standard" (i.e.
normalized) Brownian motion in $\mathbb{R}^d$.
\\
Fixing $t^*>0$ and 
$$
Y^*(a)\in \mathbb{R}^d,\;\forall a=1,\cdots,N,
$$ 
it is easy to compute the probability that, at time $t^*$, each particle
occupies the position given by $Y^*$:
$$
{\rm{Prob}}(X^\varepsilon_{t^*}(a){\approx} Y^*(a),\;\forall a=1,\cdots,N)
$$
$$
\approx
\prod_{a=1}^N [\exp(\frac{-|Y^*(a)-s_0(a)|^2}{2 \varepsilon t^*})
(2\pi \varepsilon t^*)^{-d/2}],
$$
or, in other words
$$
{\rm{Prob}}(X^\varepsilon_{t^*}{\approx} Y^*)
\approx
\exp(\frac{-||Y^*-s_0||^2}{2 \varepsilon t^*})
(2\pi \varepsilon t^*)^{-Nd/2}.
$$
Here we have denoted by $|\cdot|$ and $||\cdot||$ the Euclidean distance on
respectively  $\mathbb{R}^d$ and  $H=(\mathbb{R}^d)^N$.
We have also used symbol $\approx$ just to make notations lighter.
[What we precisely mean is:
for any Borel $\mathcal{B}$ subset of $(\mathbb{R}^d)^N$, the probability
that $X^\varepsilon_{t^*}$ belongs to $Y^*+\mathcal{B}$ is given by
$$
\int_{Y\in Y^*+\mathcal{B}}
\exp(\frac{-||Y-s_0||^2}{2 \varepsilon t^*})
(2\pi \varepsilon t^*)^{-Nd/2}dY,
$$
but we hope that our simplified notation is acceptable. Indeed, we will use
it again, without further notice.]
\subsection{First application of the large deviation principle}
We now want to compute the probability that $X^\varepsilon_{t^*}$ and $Y^*$
coincide, up to a permutation, property that we denote
by ${\rm{}}(X^\varepsilon_{t^*}){\underset{\mathrm{perm}}{\approx}}{\rm{}}(Y^*)$. 
We find
\begin{equation}\label{proba}
{\rm{Prob}}[{\rm{}}(X^\varepsilon_{t^*})
{\underset{\mathrm{perm}}{\approx}}
{\rm{}}(Y^*)]
\approx
\frac{1}{N!}\sum_{\sigma\in\mathcal{S}_N}\exp(\frac{-||Y^*\circ\sigma-s_0||^2}{2 \varepsilon t^*})
(2\pi \varepsilon t^*)^{-Nd/2}.
\end{equation}
When the level of noise $\varepsilon$ goes to zero, we immediately get 
$$
-\lim_{\varepsilon\rightarrow 0}\varepsilon\log 
{\rm{Prob}}[{\rm{}}(X^\varepsilon_{t^*}){\underset{\mathrm{perm}}{\approx}}{\rm{}}(Y^*)]
{\approx} \inf_{\sigma\in\mathcal{S}_N}
\frac{||Y^*\circ\sigma-s_0||^2}{2 t^*},
$$
which is a rather trivial illustration of the Laplace method and the large deviation principe (LDP) \cite{FW}. There is a more sophisticated aspect of the LDP: as the level of noise goes to zero,
the Brownian trajectories of the particles conditioned by
${\rm{}}(X^\epsilon_{t^*}){\underset{\mathrm{perm}}{\approx}}{\rm{}}(Y^*)$ behave more and more as
constant speed minimizing geodesic curves. More precisely, for all $t\in [0,t^*]$,
\begin{equation}\label{geod}
X_t^\varepsilon\sim_{\varepsilon\rightarrow 0}\;\; X_t=s_0+\frac{t}{t^*}(Y^*\circ \sigma^*-s_0),
\;\;
\sigma^*={\rm{Arg}}\inf_{\sigma\in\mathcal{S}_N}
\frac{||Y^*\circ\sigma-s_0||^2}{2 t^*}.
\end{equation}
which implies
\begin{equation}\label{speed}
\frac{dX_t}{dt}
=\frac{Y^*\circ \sigma^*-s_0}{t^*}
=\frac{X_{t^*}-s_0}{t^*}=\frac{X_{t}-s_0}{t},
\;\;\forall t\in ]0,t^*].
\end{equation}

Now, we are going to translate this large deviation principle into a self-consistent dynamical
system for the particles.
Let us first denote by $\pi(Y)$ (as in subsection \ref{subsection-MAG})
the unique closest point on $S$ of a "generic" point $Y\in H=(\mathbb{R}^d)^N$
\begin{equation}\label{S3}
\pi(Y)={\rm{Arg}}\inf_{s\in S}\frac{||Y-s||^2}{2}
=\nabla\Pi(Y),\;\;\;\Pi(Y)=\sup_{s\in S}\;((Y,s))-\frac{||s||^2}{2}
\end{equation}
(where $((\cdot,\cdot))$ denotes the inner product attached to $||\cdot||$ on $H$).
Then we state:
\begin{Proposition}
\label{prop}
Equation (\ref{speed}), that we have derived from (\ref{brownian})
(through a large deviation principle), can be written
as a self-consistent ordinary
differential equation for $X_t$ at least for $t>0$, 
\begin{equation}\label{GF0}
t\frac{dX_t}{dt}
=X_t-\pi(X_{t}).
\end{equation}
\end{Proposition}

\subsubsection*{Proof}
To get this result, our simple but crucial observation is that,
along the geodesic curve defined by (\ref{geod}), namely
$$ 
X_t=s_0+\frac{t}{t^*}(Y^*\circ \sigma^*-s_0),
\;\;
\sigma^*={\rm{Arg}}\inf_{\sigma\in\mathcal{S}_N}
||Y^*\circ\sigma-s_0||,
$$
$s_0$ is the closest point in $S$
not only of the end-point $X_{t^*}=Y^*\circ\sigma^*$ but also of all points
$X_t$, $\;\forall t\in [0,t^*]$.
[Then, we can write $s_0=\pi(X_t)$
in equation (\ref{speed}) which immediately leads to (\ref{GF0}) and
completes the proof of our Proposition.]
\\
Although this property is geometrically quite obvious, let us provide
a comprehensive proof for the sake of completeness.
\\
By definition (\ref{S3})
of the closest point
operator $\pi$, 
it is enough to show that
$$
\kappa=||X_t-s_0\circ\sigma||^2-||X_t-s_0||^2\ge 0
$$
for each permutation $\sigma\in\mathcal{S}_N$ and each $t\in [0,t^*]$.
By definition (\ref{geod}),
$$
\kappa=
||s_0+\frac{t}{t^*}(X^*-s_0)-s_0\circ\sigma||^2-||\frac{t}{t^*}(X^*-s_0)||^2
$$
where we set, still according to (\ref{geod}),
\begin{equation}\label{sigma}
X^*=X_{t^*}=Y^*\circ\sigma^*,
\;\;\;\sigma^*={\rm{Arg}}\inf_{\sigma\in\mathcal{S}_N}||Y^*\circ\sigma-s_0||.
\end{equation}
Expanding the squares, we obtain
$$
\kappa=
||s_0-s_0\circ\sigma||^2
+\frac{2t}{t^*}((X^*-s_0,s_0-s_0\circ\sigma)).
$$
By definition of $X^*=Y^*\circ\sigma^*$ and $\sigma^*$, we get from (\ref{sigma})
$$
||X^*-s_0||=
||Y^*\circ\sigma^*-s_0||\le ||Y^*\circ\sigma^*\circ\sigma^{-1}-s_0||
=||X^*-s_0\circ \sigma||.
$$
Thus
$
((X^*,s_0-s_0\circ\sigma))\ge 0.
$
So, we deduce
$$
\kappa\ge ||s_0-s_0\circ\sigma||^2
-\frac{2t}{t^*}((s_0,s_0-s_0\circ\sigma))
=2(1-\frac{t}{t^*})[||s_0||^2-((s_0,s_0\circ\sigma))]\ge 0
$$
(using $||s_0||=||s_0\circ\sigma||$ and the Cauchy-Schwarz inequality),
which is the desired inequality  and completes the proof of Proposition  \ref{prop}.
\\
\\
Equation (\ref{GF0}) is clearly singular at time $t=0$.
We can lift this singularity  with an exponential rescaling of 
time $t=\exp(\theta)$, $\theta\in\mathbb{R}$, and finally obtain:
\begin{equation}\label{GFLD}
\frac{dX_\theta}{d\theta}
=X_\theta-\pi(X_{\theta})=X_\theta-\nabla\Pi(X_{\theta}).
\end{equation}
So we have obtained a first order (in time) dynamical system, of "gradient type".
\subsubsection*{Miscellaneous remarks}
i) Since $\pi$ is the gradient of a Lipschitz convex function
(namely $\Pi$ defined by (\ref{Pi})), 
equation (\ref{GFLD}) is uniquely solvable in the
framework of "maximal monotone operator theory" \cite{Brz}. As already discussed in \cite{Br8},
this is a way to introduce a dissipative mechanism in the motion of particles, such as sticky collisions
when particles stick to each other while conserving their momentum (but not their kinetic
energy which decreases). At the level of the present paper, we do not want to enter such
considerations and leave (\ref{GFLD}) just as a $formal$ equation.
\\
ii) Notice that (\ref{GFLD}) is a gradient flow in the variable $X$, valued in the Hilbert space $H$, which is the counterpart of the
gradient flow in the variable Law$(X)$, in the so-called "Wasserstein space"
of half the negative squared "Wasserstein distance",
as discussed in \cite{AGS}. For us, it is important to keep a formulation in terms of
$X\in H$ and not in terms of Law$(X)$, because, in the sequel of our discussion,
it is crucial to restore the "individuality" of the 
gravitating particles during their motion,
in order to get a second order dynamical system
for them.
\\
iii) Quite remarkably, as explained in \cite{Br8}, equation (\ref{GF0}) is nothing but the
Zeldovich model used in Cosmology \cite{Ze,SZ,FMMS,BFHLMMS} as
an approximation of semi-Newtonian gravitation in an Einstein-de Sitter space!
\\
iv) Let us finally provide a possible interpretation of equation (\ref{GFLD}), vaguely
related to the so-called "anthropic principle":
What is observed at out ``present'' time $t^*$ is just a random output of 
the independent
Brownian trajectories of a large number $N$ of $indistinguishable$
particles initially located on the
set $\mathcal{A}$. As the noise vanishes,
the motion of these particles, conditioned by what we can observe now,
just looks driven by
the deterministic law (\ref{GFLD}). Of course, this is a highly speculative and questionable interpretation
coming from a mathematician and not from an authorized physicist.

\subsection{Second application of the large deviation principle}
Equation (\ref{GFLD}) is a $first$ order dynamical system, a so-called "gradient flow" since
(\ref{GFLD}) can also be written
$$
dX_\theta=\nabla\Phi(X_{\theta})d\theta,
$$
where $\Phi$ is half of the squared distance function to $S$:
\begin{equation}\label{Phi}
\Phi(Y)=\frac{||Y||^2}{2}-\Pi(Y)
=\inf_{s\in S}\frac{||Y-s||^2}{2},
\end{equation}
where $\Pi$ is defined by (\ref{Pi}). Describing a gravitational theory by a gradient flow does not
sound reasonable. We would rather like to get a $second$
order, conservative, dynamical system. Here again, the large deviation principle
turns out to be useful to get such a second order system
out of (\ref{GFLD}). Fixing $\eta>0$, we introduce the "noisy"
version of (\ref{Phi}) defined by:
\begin{equation}\label{GF-noise}
dX^\eta_\theta=\nabla\Phi(X^\eta_\theta)d\theta+\eta\; dB_\theta,
\end{equation}
where $B_\theta$ (again) denotes a Brownian process in $(\mathbb{R}^d)^N$.
Given two points $Y_0$ and $Y_1$ in $(\mathbb{R}^d)^N$, we expect from the LDP (or, more precisely,
the Freidlin-Wentzell theorem \cite{FW,DZ}) that,
as the level of noise $\eta$ goes to zero,
$$
-\eta\log 
{\rm{Prob}}[X^\eta_{\theta_0}{\underset{\mathrm{}}{\approx}}Y_0\;{\rm{and}}\;X^\eta_{\theta_1}{\underset{\mathrm{}}{\approx}}Y_1]
\sim_{\eta\rightarrow 0} \;\;\mathcal{A}(\theta_0,Y_0,\theta_1,Y_1)
$$
\begin{equation}\label{LAP}
\mathcal{A}(\theta_0,Y_0,\theta_1,Y_1)
=\inf\{\int_{\theta_0}^{\theta_1}
\frac{1}{2}||\frac{dX_\theta}{d\theta}-\nabla\Phi(X_\theta)||^2 d\theta,
\;\;X_{\theta_0}=Y_0,\;\;X_{\theta_1}=Y_1\}.
\end{equation}
In addition, as $\eta$ goes to zero,
$$
X^\eta_{\theta}\sim X_\theta,\;\;\forall\theta\in [\theta_0,\theta_1],
$$
where
$$
X={\rm{Arg}}\inf\{\int_{\theta_0}^{\theta_1}
\frac{1}{2}||\frac{dX_\theta}{d\theta}-\nabla\Phi(X_\theta)||^2 d\theta,
\;\;X_{\theta_0}=Y_0,\;\;X_{\theta_1}=Y_1\}.
$$
Strictly speaking, this is correct when $\nabla\Phi$ is Lipschitz continuous,
which is not true in our case (where $\nabla\Phi$ is not even continuous),
without further assumptions on the data.
However, from the pure modeling viewpoint, it is very tempting to find the second
order dynamical system linked to the least action principle (\ref{LAP}).
Since $\Phi$, as defined by (\ref{Phi}), is half of a squared distance function, we have
\begin{equation}\label{HJ}
\frac{1}{2}||\nabla\Phi(Y)||^2=\Phi(Y),
\end{equation}
for every $Y\in H\setminus \mathcal{N}$, where $\mathcal{N}$ is the set on which $\Phi$
is not differentiable, which is a negligible subset
of $H$, both in the Lebesgue almost everywhere sense
and in the Baire category sense.
Thus, at least for each curve $X$ that stays away from $\mathcal{N}$ for Lebesgue
almost every $\theta\in [\theta_0,\theta_1]$, we have
$$
\int_{\theta_0}^{\theta_1}
\frac{1}{2}||\frac{dX_\theta}{d\theta}-\nabla\Phi(X_\theta)||^2 d\theta
$$
$$
=\int_{\theta_0}^{\theta_1}
(\frac{1}{2}||\frac{dX_\theta}{d\theta}||^2+||\nabla\Phi(X_\theta)||^2)d\theta
-\int_{\theta_0}^{\theta_1}
\nabla\Phi(X_\theta)\cdot\frac{dX_\theta}{d\theta}d\theta
$$
$$
=\int_{\theta_0}^{\theta_1}
(\frac{1}{2}||\frac{dX_\theta}{d\theta}||^2+\Phi(X_\theta))d\theta
-\Phi(X_{\theta_1})+\Phi(X_{\theta_0}).
$$
Therefore, the action principle (\ref{LAP}) is equivalent to
\begin{equation}\label{LAP2}
\mathcal{\tilde A}(\theta_0,Y_0,\theta_1,Y_1)
=\inf\{\int_{\theta_0}^{\theta_1}
(\frac{1}{2}||\frac{dX_\theta}{d\theta}||^2+\Phi(X_\theta)) d\theta,
\;\;X_{\theta_0}=Y_0,\;\;X_{\theta_1}=Y_1\}.
\end{equation}
From this equivalent LAP, we find as optimality equation
\begin{equation}\label{dynamical}
\frac{d^2 X_\theta}{d\theta^2}=\nabla\Phi(X_\theta)=X_\theta-\pi(X_{\theta}),
\end{equation}
which is just the second order version of the gradient flow equation (\ref{GFLD}).
More explicitly, we have obtained
\begin{equation}\label{dynamical2}
\frac{d^2 X_\theta(a)}{d\theta^2}=X_\theta(a)-A(\sigma_\theta(a)),
\end{equation}
where 
$$
\sigma_\theta={\rm{Arg}}\inf_{\sigma\in \mathcal{S}_N}
\sum_{a=1}^N\;|X_\theta(a)-A(\sigma(a))|^2,
$$
which is nothing but the discrete version (\ref{MAG3}) of 
the Monge-Amp\`ere model of
gravitation discussed in the first section. 
\\
So we have achieved, at least
at a formal level, the main goal of our paper which was,
through a double application of the large deviation principle,
the derivation of the (discrete)
Monge-Amp\`ere gravitational model 
from one of the simplest thinkable model of particles:
$N$ independent Brownian trajectories, with
an intriguing interplay between their indistinguishability
(for the first application of the LDP) and their individuality (for
the second application of the LDP).
\subsubsection*{Acknowledgments}
This work has been partly supported by the grant ISOTACE
ANR-12-MONU-0013 (2012-2016).
The author wants to thank the hospitality of the Radon Institute (RICAM, Linz) and the
Hausdorff Institute (HIM, Bonn) during
the special semester on New Trends in Calculus of Variations (Oct-Dec 2014) and
the Junior Hausdorff trimester program on Optimal Transportation
(Jan-Apr 2015). 
\\
http://www.ricam.oeaw.ac.at/specsem/specsem2014
\\
http://wt.iam.uni-bonn.de/conference-new-trends-in-optimal-transport/
\\
He also thanks Christian L\'eonard for very stimulating
discussions and comments bringing together large deviation and optimal transport theories.

\end{document}